\documentclass{article}       % onecolumn (second format)
\usepackage{anysize}
\marginsize{2.5cm}{2.5cm}{2cm}{2cm}
\usepackage{graphicx}
\usepackage{amssymb}
\usepackage{makeidx}         % allows index generation
\usepackage{amsmath}
\usepackage{multicol}        % used for the two-column index
\usepackage[bottom]{footmisc}% places footnotes at page bottom
\usepackage{color}
\usepackage{amsfonts}
\usepackage[latin1]{inputenc}
\usepackage{epstopdf}
\usepackage{url}

\newcommand{\eq}{\begin{equation}}
\newcommand{\eeq}{\end{equation}}
\newcommand{\eqn}{\begin{eqnarray}}
\newcommand{\eeqn}{\end{eqnarray}}

\newcommand{\bsea}{\begin{subeqnarray}}
\newcommand{\esea}{\end{subeqnarray}}
\newcommand{\nn}{\nonumber}

\newcommand{\Sp}[2]{\left< #1,#2 \right> }

\newcommand{\Min}[1]{\,\underset{#1}{\mathrm{min}}\,}

\newcommand{\Lim}[1]{\,\underset{#1}{\mathrm{lim}}\,}

\newcommand{\de}{\mathrm{d}}
\newcommand{\tr}{\mathop{\rm tr}}  %traccia
\newcommand{\e}[1]{\mathrm{e}^{#1}}
\newcommand{\Set}[1]{\left\{ #1\right\}}

\newcommand{\Bc}{ \mathcal{B}}
\newcommand{\Cc}{ \mathcal{C}}

\newcommand{\Hc}{ \mathcal{H}}
\newcommand{\Ic}{ \mathcal{I}}

\newcommand{\Kc}{ \mathcal{K}}
\newcommand{\Lc}{ \mathcal{L}}

\newcommand{\Qc}{ \mathcal{Q}}

\newcommand{\Sc}{ \mathcal{S}}

\newcommand{\Vc}{ \mathcal{V}}

\newcommand{\Es}{ \mathbb{E}}

\newcommand{\Ns}{ \mathbb{N}}

\newcommand{\Rs}{ \mathbb{R}}
\newcommand{\Ss}{ \mathbb{S}}
\newcommand{\Ts}{ \mathbb{T}}

\newcommand{\Zs}{ \mathbb{Z}}

\def\qed{\hfill \vrule height 7pt width 7pt depth 0pt \smallskip}

\newtheorem{remark}{Remark}[section]
\newtheorem{theorem}{Theorem}[section]

\newtheorem{proposition}{Proposition}[section]
\newtheorem{lemma}{Lemma}[section]
\newtheorem{problem}{problem}[section]
\newcommand{\proof}{\noindent {\bf Proof. }}

\newcommand{\teo}{\begin{theorem}}
\newcommand{\eteo}{\end{theorem}}
\newcommand{\cor}{\begin{corollary}}
\newcommand{\ecor}{\end{corollary}}
\newcommand{\prop}{\begin{proposition}}
\newcommand{\eprop}{\end{proposition}}
\newcommand{\lem}{\begin{lemma}}
\newcommand{\elem}{\end{lemma}}
\newcommand{\pb}{\begin{problem}}
\newcommand{\epb}{\end{problem}}
\newcommand{\rem}{\begin{remark}}
\newcommand{\erem}{\end{remark}}
\newcommand{\ej}{\mathrm{e}^{j\vartheta}}

\newcommand{\Sts}{\Ss_+(\Ts)}

\newcommand{\Lm}{\mathrm{L}_\infty(\Ts)}
\newcommand{\Dkl}{\Sc_{\mathrm{KL}}}

\newcommand{\Da}{\Sc_\alpha}

\newcommand{\Danu}{\Sc_{\nu}}

\newcommand{\Rgamma}{\mathrm{Range}\;\Gamma}

\newcommand{\Lag}{\Lc^\Gamma}
\newcommand{\Lbg}{\Lc_{\nu}^\Gamma}

\begin{document}

\title{Rational approximations of spectral densities\\ based on the Alpha divergence}

\author{Mattia Zorzi\thanks{M. Zorzi
is with the Department of Electrical Engineering and Computer
Science, University of Liège, 4000 Liège, Belgium (e-mail:
mzorzi@ulg.ac.be)}}
\date{}

\maketitle

\begin{abstract}
We approximate a given rational spectral density by one that is
consistent with prescribed second-order statistics. Such an
approximation is obtained by selecting the spectral density having
minimum ``distance'' from under the constraint corresponding to
imposing the given second-order statistics. We analyze the
structure of the optimal solutions as the minimized ``distance''
varies in the Alpha divergence family. We show that the
corresponding approximation problem leads to a family of rational
solutions. Secondly, such a family contains the solution which
generalizes the Kullback-Leibler solution proposed by Georgiou and
Lindquist in 2003. Finally, numerical simulations suggest that
this family contains solutions close to the non-rational solution
given by the principle of minimum discrimination information.
\end{abstract}

\section{Introduction}  \label{section_THREE_est}
This paper deals with the rational approximation of power spectra
of stationary stochastic processes. More precisely, we consider
the following situation. Let $y=\Set{y_k;\; k\in\Zs}$ be a zero
mean, $\Rs$-valued, purely non-deterministic, full-rank,
stationary process with unknown spectral density $\Omega(\ej)$
defined on the unit circle $\Ts$. A rational {\em prior} spectral
density $\Psi\in\Sts$, which represents the {\em a priori}
information on $y$, is available. Here, $\Sts$ denotes the family
of bounded and coercive $\Rs$-valued spectral density functions on
$\Ts$. Then, some second-order statistics of $y$ are observed.
These are encoded in the output covariance, denoted by $\Sigma$,
of a rational filters bank $G(z)=(zI-A)^{-1}B$ driven by $y$. The
filter parameters are chosen in such a way that $A\in\Rs^{n\times
n}$ is a stability matrix, $B\in\Rs^{n}$, and $(A,B)$ is a
reachable pair with $n>1$. Thus the output process, denoted by
$x=\{x_k;\; k\in\Zs\}$, is a stationary process, and
$\Sigma=\Es[x_kx_k^T]=\int G\Omega G^*$ is a positive definite
matrix. Here, integration takes places on $[-\pi,\pi)$ with
respect to the normalized {\em Lebesgue} measure $\de\vartheta/
2\pi$. The rational prior $\Psi$ is typically inconsistent with
$\Sigma$, i.e. $\Sigma\neq \int G\Psi G^*$. Hence, our task
consists in finding a rational spectral density $\Phi\in\Sts$
which is as close as possible to $\Psi$ and such that
\eq\label{vincolo_Sigma}\Sigma=\int G\Phi G^*.\eeq The closeness
between $\Phi$ and $\Psi$ is quantified by considering a distance
measure $\Sc(\cdot\|\cdot)$ among spectral densities in
$\Ss_+(\Ts)$, i.e. $\Sc(\Phi_1\|\Phi_2)\geq 0$ for each
$\Phi_1,\Phi_2\in\Ss_+(\Ts)$ and equality holds if and only if
$\Phi_1=\Phi_2$ but we do not require that symmetry and triangular
inequality hold. In other words we reconcile the inconsistency
among $\Sigma$ and $\Psi$ by approximating the spectrum $\Psi$ by
a rational spectrum $\Phi$ compatible with $\Sigma$. This problem
is motivated by THREE-like spectral estimation paradigms
\cite{A_NEW_APPROACH_BYRNES_2000,KL_APPROX_GEORGIUO_LINDQUIST,PAVON_FERRANTE_ONGEORGIOULINDQUIST,Hellinger_Ferrante_Pavon,MATRICIAL_ALGTHM_RAMPONI_FERRANTE_PAVON_2009,FERRANTE_ON_THE_CONVERGENCE,FERRANTE_TIME_AND_SPECTRAL_2012,BETA}
wherein $\Phi$ represents an estimate of the unknown spectral
density $\Omega$ compatible with the given second-order statistics
of $\Omega$ and as close as possible to the {\em
 a priori} information given.

The role of the filters bank $G(z)$ consists in providing the
interpolation conditions for the solution to the spectrum
approximation problem. More specifically, by choosing the filters
bank poles appropriately, it is possible to give preference to
selected frequency bands of $\Phi$ and allow more accurate
reconstruction of $\Omega$ in these particular frequency bands,
\cite{A_NEW_APPROACH_BYRNES_2000}.

Concerning the distance measures in the above problem, we have
many options. In \cite{KL_APPROX_GEORGIUO_LINDQUIST}, the {\em
Kullback-Leibler} divergence among spectral densities having the
same zeroth moment is considered \eq\label{KL_0}
\Sc_{\mathrm{KL0}}(\Phi_1\|\Phi_2):=\int
\Phi_1\log\left(\frac{\Phi_1}{\Phi_2}\right).\eeq The found
rational solution $\Phi$ is the spectral density satisfying
(\ref{vincolo_Sigma}) which minimizes
$\Sc_{\mathrm{KL0}}(\Psi\|\Phi)$. Note that in statistics,
information theory and communications, the {\em Kullback-Leibler}
divergence is typically minimized with respect to the first
argument, \cite{COVER_THOMAS,CSISZAR_MATUS}, according to the {\em
Principle of minimum discrimination information}, MinxEnt,
\cite{KULLBACK59}. It states that, given new information, a new
distribution $\Phi$ should be chosen in such a way as to minimize
$\Sc_{\mathrm{KL0}}(\Phi\|\Psi)$. In
\cite{KL_APPROX_GEORGIUO_LINDQUIST}, the unusual choice to
minimize with respect to the second argument is dictated by the
fact that the optimization of $\Sc_{\mathrm{KL0}}(\Phi\|\Psi)$
leads to a non-rational ``exponential-type'' solution,
\cite{RELATIVE_ENTROPY_GEORGIOU_2006}. Then, in
\cite{Hellinger_Ferrante_Pavon} and
\cite{FERRANTE_TIME_AND_SPECTRAL_2012}, the rational solutions
corresponding to the {\em Hellinger} distance and the {\em
Itakura-Saito} distance, respectively, are considered. Finally, in
\cite{BETA} the Beta divergence family is considered. The latter
smoothly connects the {\em Kullback-Leibler} divergence with the
{\em Itakura-Saito} distance. Making additional assumptions on
$\Psi$ besides the rationality, it is possible to prove that the
Beta divergence leads to a family of rational solutions. Finally,
it is worth noting that the solutions considered in
\cite{Hellinger_Ferrante_Pavon,FERRANTE_TIME_AND_SPECTRAL_2012,BETA}
also hold for the multichannel case, i.e. $y$ is a multivariate
process.

 The aim of this paper is to solve the previous spectrum approximation
 problem by employing the Alpha divergence family which
 smoothly connects the {\em Kullback-Leibler} divergence
\eq \Dkl(\Phi_1\|\Phi_2):=\int\Phi_1
\log\left(\frac{\Phi_1}{\Phi_2}\right)-\Phi_1+\Phi_2\eeq with its
 reverse $\Sc_{\mathrm{KL}}(\Phi_2\|\Phi_1)$, passing through the {\em Hellinger} distance.
Firstly we will generalize the solution in
\cite{KL_APPROX_GEORGIUO_LINDQUIST} to spectral densities with
different zeroth moment. Then, we will see that it is possible to
characterize a family of rational solutions to the problem without
making additional assumptions on $\Psi$. The limit of this family
of solutions has the same ``exponential-type'' structure as the
MinxEnt solution. Furthermore, simulations suggest that the limit
above converges to the non-rational solution given by
$\Sc_{\mathrm{KL}}(\Phi\|\Psi)$. Thus, the spectrum approximation
problem based on the Alpha divergence family also seems to provide
rational solutions, for a suitable choice of the $\alpha$
parameter, close to the MinxEnt solution.

\section{Spectrum approximation problem and feasibility conditions}\label{section_feasibility}
The spectrum approximation problem we are dealing with minimizes a
suitable distance measure, denoted by $\Sc(\Phi\|\Psi)$, over the
set \eq \Ic_\Sigma:=\left\{\Phi\in\Sts\;|\; \int G\Phi
G^*=\Sigma\right\}\eeq where $\Psi\in\Sts$ is rational, and
$G(z)=(zI-A)^{-1}B$ has the same properties specified in the
Introduction. Both $\Psi$ and $G(z)$ are given. The output
covariance $\Sigma$ is estimated by a given finite-length sequence
$y_1\ldots y_N$, extracted from a realization of $y$, see
\cite{ME_ENHANCEMENT_FERRANTE_2012,ON_THE_ESTIMATION_ZORZI_2012}.
Thus, the first problem concerns the feasibility of the spectrum
approximation problem, i.e. when the set $\Ic_\Sigma$ is non-empty
for a given $\Sigma$. To deal with this issue, we first introduce
some notation: $\Qc_{n}\subset \Rs^{n\times n}$ denotes the
$n(n+1)/2$-dimensional real vector space of $n$-dimensional
symmetric matrices and $\Qc_{n,+}$ denotes the corresponding cone
of positive definite matrices. We denote as $\Vc(\Sts)$ the linear
space generated by $\Sts$. Finally, we introduce the linear
operator
\eqn \Gamma &:& \Vc(\Sts)\rightarrow \Qc_{n} \nn\\
&& \Phi\mapsto \int G\Phi G^*.\eeqn In
\cite{GEORGIUO_THE_STRUCUTRE_2002}, it was shown that a matrix
$P\in\Qc_{n}$ belongs to the range of $\Gamma$, denoted by
$\Rgamma$, if and only if there exists $H\in \Rs^{1\times
    n}$ such that \eq \label{equazione_R_gamma}P-A P A^T=BH+H^TB^T .\eeq
It turns out that the spectrum approximation problem is feasible
if and only if $\Sigma\in\Rgamma\cap \Qc_{n,+}$, see
\cite{GEORGIUO_THE_STRUCUTRE_2002,Hellinger_Ferrante_Pavon}.
 Once we have $\Sigma$ in
such a way that the spectrum approximation problem is feasible, we
can replace $G$ with $\overline{G}=\Sigma^{-\frac{1}{2}}G$ and
$(A,B)$ with
$(\overline{A}=\Sigma^{-\frac{1}{2}}A\Sigma^{\frac{1}{2}},\overline{B}=\Sigma^{-\frac{1}{2}}B)$.
Thus, constraint (\ref{vincolo_Sigma}) may be rewritten as
\eq\label{vincolo_I}\int \bar{G}\Phi \bar{G}^*=I.\eeq Accordingly,
from now on we assume that our spectrum approximation problem is
feasible and we consider the following equivalent formulation:
Given a rational {\em prior} $\Psi\in\Sts$ and $G(z)$ such that
$I\in\Rgamma$, minimize $\Sc(\Phi\|\Psi)$ over the (non-empty) set
\eq \left\{\Phi\in\Sts\;|\;\int G\Phi G^*=I\right\}.\eeq In the
following section we will introduce the last element
characterizing our problem: The distance measure
$\Sc(\Phi\|\Psi)$.

\section{Spectrum approximation problem with the Alpha divergence}\label{section_THREE_alpha}
The Alpha divergence family,
\cite{AMARI_DIFFERENTIAL_GEOMETRIC,FLEXIBLE_ROBUST_CICHOKI_AMARI_2010},
among two spectral densities $\Phi_1,\Phi_2\in\Sts$ is defined as
\eqn \Da(\Phi_1\|\Phi_2):= \int
\frac{1}{\alpha(\alpha-1)}\Phi_1^{\alpha}\Phi_2^{1-\alpha}-\frac{1}{\alpha-1}\Phi_1+\frac{1}{\alpha}\Phi_2\eeqn
where $\alpha$ is a real parameter. For $\alpha=0$ and $\alpha=1$,
it is defined by continuity \eq \lim_{\alpha\rightarrow
0}\Da(\Phi_1\|\Phi_2)=\Dkl(\Phi_2\|\Phi_1),\;\;\lim_{\alpha\rightarrow
1}\Da(\Phi_1\|\Phi_2)=\Dkl(\Phi_1\|\Phi_2). \eeq Thus, the Alpha
divergence is a continuous function of real variable $\alpha$ in
the whole range including singularities and it smoothly connects
$\Dkl(\Phi_1\|\Phi_2)$ with its reverse $\Dkl(\Phi_2\|\Phi_1)$.
Moreover, $\Sc_\alpha$ is strictly convex with respect to both
$\Phi_1$ and $\Phi_2$. Note that, for $\alpha=\frac{1}{2}$ we
obtain, up to a constant factor, the {\em Hellinger} distance \eq
\Sc_{\mathrm{H}}(\Phi_1,\Phi_2):=\int
(\sqrt{\Phi_1}-\sqrt{\Phi_2})^2\eeq and for $\alpha=2$ we obtain
the {\em Pearson Chi-square} distance \eq
\Sc_{\mathrm{P}}(\Phi_1\|\Phi_2):=\frac{1}{2}\int\frac{(\Phi_1-\Phi_2)^2}{\Phi_2}.\eeq

Since $\alpha\in\Rs$, we choose the parametrization
$\alpha=-\frac{1}{\nu}+1$ with $\nu\in\Zs\setminus\Set{0}$ and we
consider the following spectrum approximation problem. \pb
\label{problema_THREE_I} Given a rational {\em prior}
$\Psi\in\Sts$ and $G(z)=(zI-A)^{-1}B$ such that $I\in\Rgamma$,
\eqn \label{constraint_pb_THREE_I}&&\mathrm{minimize}\; \;
\Danu(\Phi\|\Psi)\; \; \mathrm{over \; the\; set}\nn\\ &&
\hspace{2cm}\; \; \Set{\Phi\in\Sts\;|\; \int G\Phi G^*=I},\eeqn
where
 \eq \Danu(\Phi\|\Psi):=\left\{
                                    \begin{array}{ll}
                                    \Sc_{\mathrm{KL}}(\Psi\|\Phi) ,& \nu=1 \\
                                    \int
\frac{\nu^2}{1-\nu}\Phi^{\frac{\nu-1}{\nu}}\Psi^{\frac{1}{\nu}}+\nu\Phi+\frac{\nu}{\nu-1}\Psi
  \hspace{0.5cm}& \nu\in \Zs\setminus\Set{0,1}\\
\Sc_{\mathrm{KL}}(\Phi\|\Psi)  &\nu= \infty
                                    \end{array}
                                  \right.
\eeq \epb The above problem is a constrained convex optimization
problem and admits at most one solution because
$\Sc_\nu(\cdot\|\Psi)$ is strictly convex over $\Ss_+(\Ts)$. On
the other hand the existence issue is not trivial. In fact, since
the set $\Set{\Phi\in\Sts\;|\; \int G\Phi G^*=I}$ is open, it is
possible that the minimum point of $\Sc_\nu$ is not attained. In
the next sections we show that Problem \ref{problema_THREE_I}
admits a unique solution once fixed $\nu\in\Zs$ such that $1 \leq
\nu\leq  \infty$. This task is accomplished by exploiting the
duality theory, \cite{BOYD_CONVEX_OPTIMIZATION}. We not only show
that the dual problem admits solution but we even prove that there
exists a suitable subset of the {\em dual functional} domain which
contains a unique optimal Lagrange multiplier for Problem
\ref{problema_THREE_I}. Accordingly, it is possible to employ the
efficient Newton-type algorithm presented in
\cite{MATRICIAL_ALGTHM_RAMPONI_FERRANTE_PAVON_2009} for computing
the numerical solution of the problem. Moreover, thanks to the
chosen parametrization $\nu$, the duality theory will show that
the solution is rational with a prescribed maximum degree when $1
\leq \nu<\infty$. Thus, Problem \ref{problema_THREE_I} with $1\leq
\nu<\infty$ can be viewed as a {\em Nevanlinna-Pick} interpolation
problem with bounded degree,
\cite{Blomqvist_MATRIX_VALUED_2003,INTERPOLATION_GEORGIUO_1999}.
When $\nu\in\Zs$ is such that $\nu\leq -1$, however, the existence
of the optimal solution is not guaranteed, see Remark
\ref{remark_nu_neg}.

The analysis will be divided in the following three cases:
$\nu=1$, $1<\nu<\infty$ and $\nu=\infty$.

\section{Case $\nu=1$}\label{section_THREE_KL1} The corresponding
{\em Lagrangian} is \eqn
L(\Phi, \Lambda)&:=&\Sc_{\mathrm{KL}}(\Psi\|\Phi)+\int \Psi+\Sp{\int G\Phi G^*-I}{\Lambda}\nn\\
&=&\int \left[-\Psi\log(\Psi^{-1}\Phi)+\Phi+G^*\Lambda
G\Phi\right]-\tr(\Lambda)\eeqn where we exploited the fact that
$\int\Psi$ plays no role in the optimization task. The Lagrange
multiplier $\Lambda\in\Qc_n$ can be uniquely decomposed as
$\Lambda=\Lambda_\Gamma+\Lambda_\bot$ where
$\Lambda_\Gamma\in\Rgamma$ and $\Lambda_\bot\in[\Rgamma]^\bot$.
Since $\Lambda_\bot$ is such that $G^*(\ej)\Lambda_\bot
G(\ej)\equiv 0$ and $\tr(\Lambda_\bot)=\Sp{\Lambda_\bot}{I}=0$,
see \cite[Section
III]{MATRICIAL_ALGTHM_RAMPONI_FERRANTE_PAVON_2009}, it does not
affect the {\em Lagrangian}, i.e.
$L(\Phi,\Lambda)=L(\Phi,\Lambda_\Gamma)$. Accordingly we can
impose from now on that $\Lambda\in\Rgamma$.

Consider now the unconstrained minimization problem
\eq\min_{\Phi}\Set{L(\Phi,\Lambda)\;|\; \Phi\in\Sts}.\eeq Since
$\Sc_{\mathrm{KL}}(\Psi\|\Phi)$ is strictly convex with respect to
$\Phi\in\Sts$, then $L(\cdot,\Lambda)$ is strictly convex with
respect to $\Phi\in\Sts$. Accordingly, the unique minimum point of
$L(\cdot,\Lambda)$ is given by annihilating its first variation in
each direction $\delta\Phi\in L_\infty(\Ts)$: \eqn
\label{variazione_delta_phi_KL}\delta L(\Phi,\Lambda;\delta
\Phi)&=&\int\left[\left(-\Psi\Phi^{-1}+1+ G^*\Lambda
G\right)\delta\Phi\right]. \eeqn Note that, $-\Psi\Phi^{-1}+1
+G^*\Lambda G\in L_\infty(\Ts)$. Hence,
(\ref{variazione_delta_phi_KL}) is zero $\forall \delta\Phi\in
L_\infty(\Ts)$ if and only if \eq -\Psi\Phi^{-1}+1+ G^*\Lambda
G\equiv 0.\eeq Therefore, the unique minimum point of
$L(\cdot,\Lambda)$ has the form \eq\Phi(\Lambda):=\frac{\Psi}{1+
G^*\Lambda G}.\eeq As $\Psi\Phi^{-1}\in\Sts $, the set of
admissible Lagrange multipliers  is
\eq\Lag:=\Set{\Lambda\in\Qc_n\;|\; 1+G^*\Lambda G>0\; \hbox{on}
\;\Ts}\cap\Rgamma.\eeq Since $\Phi(\Lambda)$ is the unique minimum
point of $L(\cdot,\Lambda)$, we get \eq
\label{disug_lagrangiane_KL}L(\Phi(\Lambda),\Lambda)<
L(\Phi,\Lambda),\; \; \forall \; \Phi\in\Sts \hbox{ s.t. }\Phi\neq
\Phi(\Lambda).\eeq Hence, if we produce $\Lambda^\circ\in\Lag$
satisfying constraint in (\ref{constraint_pb_THREE_I}), inequality
(\ref{disug_lagrangiane_KL}) implies \eq
\Sc_{\mathrm{KL}}(\Psi\|\Phi(\Lambda^\circ))\leq\Sc_{\mathrm{KL}}(\Psi\|\Phi),\;
\; \forall  \Phi\in\Sts \;\hbox{s.t.}\; \int G\Phi G^*=\Sigma\eeq
and equality holds if and only if $\Phi= \Phi(\Lambda^\circ)$.
Accordingly, such a $\Phi(\Lambda^\circ)$ is the unique solution
to Problem \ref{problema_THREE_I}. Note that, $\Phi(\Lambda)$ is
rational because $\Psi(z)$ is a rational function and $G(z)$ is a
rational matrix function. Furthermore, it is possible to
characterize an upper bound on its degree:
\eq\deg[\Phi(\Lambda)]\leq \deg[\Psi]+2n.\eeq The following step
consists in showing the existence of such a $\Lambda^\circ$ by
exploiting the duality theory. The dual problem consists in
maximizing the functional
 \eq \inf_{\Phi}
L(\Phi,\Lambda)=L(\Phi(\Lambda),\Lambda) =\int\Psi\log (1+
G^*\Lambda G)+\Psi-\tr(\Lambda)\eeq which is equivalent to
minimize the following functional hereafter referred to as {\em
dual functional}: \eq \label{funz_dual_KL}J(\Lambda)=\int-\Psi\log
(1+G^*\Lambda G)+\tr(\Lambda).\eeq \teo
\label{teorema_J_stricly_convex_KL} The dual functional $J$
belongs to $\Cc^2(\Lag)$ and it is strictly convex over
$\Lag$.\eteo \proof The first variation of $J(\Lambda)$ in
direction $\delta \Lambda_1 \in\Qc_n$ is \eq
\label{gradiente_J_KL} \delta J(\Lambda;\delta\Lambda_1)=-\int
\Psi(1+ G^*\Lambda G)^{-1}G^*\delta\Lambda_1 G+\tr(\delta
\Lambda_1).\eeq The linear form $\nabla J_{\Lambda}(\cdot):=
\delta J(\Lambda;\cdot)$ is the {\em gradient} of $J$ at
$\Lambda$. In order to prove that $J(\Lambda)\in\Cc^1(\Lbg)$ we
have to show that $\delta J(\Lambda;\delta\Lambda_1)$, for any
fixed $\delta \Lambda_1$, is continuous in $\Lambda$. To this aim,
consider a sequence $M_n\in\Rgamma$ such that $M_n\rightarrow 0$
and define $Q_N(z)=(1+G(z)^*NG(z))^{-1}$ with $N\in\Qc_n$. By
Lemma 5.2 in \cite{MATRICIAL_ALGTHM_RAMPONI_FERRANTE_PAVON_2009},
$Q_{\Lambda+M_n}$ converges uniformly to $Q_{\Lambda}$. Thus,
applying the bounded convergence theorem, we obtain \eq
\lim_{n\rightarrow \infty } \int G \Psi Q_{\Lambda+M_n} G^*= \int
G \Psi Q_{\Lambda} G^*. \eeq Accordingly,
 $\delta J(\Lambda;\delta\Lambda)$ is continuous in $\Lambda$, i.e. $J$ belongs to $\Cc^1(\Lag)$.
  The second variation in direction
$\delta\Lambda_1,\delta\Lambda_2\in\Qc_{n}$ is \eq
\label{hessiano_J_KL}\delta^2
J(\Lambda;\delta\Lambda_1,\delta\Lambda_2)=\int \Psi(1+ G^*\Lambda
G)^{-2}G^*\delta\Lambda_1 GG^*\delta\Lambda_2 G.\eeq The bilinear
form $\Hc_{\Lambda}(\cdot,\cdot)=\delta^2 J(\Lambda;\cdot,\cdot)$
is the {\em Hessian} of $J$ at $\Lambda$. The continuity of
$\delta^2 J$ can be established by using the previous
argumentation. Finally, it remains to be shown that $J$ is
strictly convex on the open set $\Lag$. Since $J\in\Cc^2(\Lag)$,
it is sufficient to show that $\Hc_{\Lambda}(\delta
\Lambda,\delta\Lambda)\geq 0$ for each $\delta \Lambda\in\Rgamma$
and equality holds if and only if $\delta\Lambda=0$. Since the
integrand in (\ref{hessiano_J_KL}) is a nonnegative function when
$\delta\Lambda_1=\delta\Lambda_2$, we have
$\Hc_{\Lambda}(\delta\Lambda,\delta\Lambda)\geq 0$. If
$\Hc_{\Lambda}(\delta\Lambda,\delta\Lambda)=0$, then
$G^*\delta\Lambda G\equiv 0$ namely $\delta
\Lambda\in[\Rgamma]^\bot$, see \cite[Section
III]{MATRICIAL_ALGTHM_RAMPONI_FERRANTE_PAVON_2009}. Since
$\delta\Lambda\in\Rgamma$, it follows that $\delta\Lambda =0$.
 In conclusion,
the Hessian is positive definite and the dual functional is
strictly convex on $\Lag$. \qed\\
In view of Theorem \ref{teorema_J_stricly_convex_KL}, the dual
problem $\Min{\Lambda}\Set{J(\Lambda)\;|\; \Lambda\in\Lag}$ admits
at most one solution $\Lambda^\circ$. Since $\Lag$ is an open set,
such a $\Lambda^\circ$ (if it does exist) annihilates the first
directional derivative (\ref{gradiente_J_KL}) for each
$\delta\Lambda\in\Qc_n$
\eq\Sp{I-\int[G\frac{\Psi}{(1+G^*\Lambda^\circ
G)}G^*}{\delta\Lambda}=0,\;\; \forall \delta\Lambda\in\Qc_n\eeq
or, equivalently, \eq \label{condizione_der_prima_kl_raz}I=\int
G\frac{\Psi}{(1+G^*\Lambda^\circ G)}G^*=\int
G\Phi(\Lambda^\circ)G^*.\eeq This means that
$\Phi(\Lambda^\circ)\in\Sts$ satisfies constraint in
(\ref{constraint_pb_THREE_I}) and $\Phi(\Lambda^\circ)$ is
therefore the unique solution to Problem \ref{problema_THREE_I}.

Although the existence question is quite delicate, since set
$\Lag$ is open and unbounded, we now show that such a
$\Lambda^\circ$ minimizing $J$ over $\Lag$ does exist. \teo The
dual functional $J$ has a unique minimum point in $\Lag$.\eteo
\proof In view of Theorem \ref{teorema_J_stricly_convex_KL} it is
sufficient to show that $J$ takes a minimum value over $\Lag$.
Consider the closure of
$\Lag$\eq\overline{\Lag}=\{\Lambda\in\Rgamma\;|\; 1+G^*\Lambda
G\geq 0\; \hbox{on} \;\Ts\}\eeq and define the sequence of
functions on $\overline{\Lag}$\eq J^n(\Lambda):=\int -\Psi
\log\left(1+G^*\Lambda G+\frac{1}{n}\right)+\tr(\Lambda).\eeq By
Lemma 1, Lemma 2 and Lemma 3 in
\cite{FERRANTE_PAVON_RAMPONI_FURTHERRESULTS} we conclude that:
\begin{itemize}
\item the pointwise limit $J^\infty(\Lambda)=\lim_{n\rightarrow
\infty}J^n(\Lambda)$ exists and is a lower semi-continuous, convex
function on $\overline{\Lag}$ with values in the extended reals
    \item $J^\infty$ is bounded below on $\overline{\Lag}$
    \item $J^\infty(\Lambda)=J(\Lambda)$ on $\Lag$
    \item $J^\infty(\Lambda)$ is finite on $\mathcal{B}^c$ which is the complement set of
    $\Bc:=\{\Lambda\in\partial\Lag\;|\; 1+G^*\Lambda G\equiv0\}$
    \item $\lim_{\|\Lambda\|\rightarrow\infty} J(\Lambda)=\infty.$
\end{itemize}
Thus, $J$ is {\em inf-compact} over $\overline{\Lag}$ and it
admits a minimum point in $\Lambda^\circ\in \overline{\Lag}$.
Clearly, $J(\Lambda)=\infty$ for $\Lambda\in\Bc$. Let
$\overline{\Lambda}\in\Bc^c$, thus $J(\overline{\Lambda})$ is
finite. Since $I\in\Lag$, then
$\overline{\Lambda}+\varepsilon(I-\overline{\Lambda})\in\overline{\Lag}$
with $\varepsilon\in[0,1]$. The one-side directional derivative is
\eq\lim_{\varepsilon \searrow 0}
\frac{J(\overline{\Lambda}+\varepsilon(I-\overline{\Lambda}))-J(\overline{\Lambda})}{\varepsilon}=
\int\Psi-\int\Psi\frac{G^*G+1}{1+G^*\overline{\Lambda}G}+\tr(I-\overline{\Lambda})=-\infty.\eeq
Thus, $\overline{\Lambda}$ cannot be a minimum point. We conclude
that $\Lambda^\circ\in\Lag$.
\qed\\
Now, we analyze the differences among our solution with $\nu=1$
and the solution given in \cite{KL_APPROX_GEORGIUO_LINDQUIST}. The
latter is obtained by minimizing $\Sc_{\mathrm{KL0}}(\Psi\|\Phi)$
in (\ref{KL_0}), and the corresponding optimal form is \eq
\Phi_{\mathrm{KL0}}(\Lambda)=\frac{\Psi}{G^*\Lambda G}.\eeq As
noticed in the Introduction, $\Sc_{\mathrm{KL0}}(\Psi\|\Phi)$ is a
distance measure among spectral densities having the same zeroth
moment. However, when the matrix $A$ is singular the zeroth moment
of $\Phi$ is fixed by constraint
(\ref{constraint_pb_THREE_I}):\eq\int
\Phi=\frac{v^*v}{\|v^*B\|^2},\eeq where $v\in\Rs^n$ is such that
$v^*A=0$, see \cite{KL_APPROX_GEORGIUO_LINDQUIST}. Hence, \eq
\Sc_{\mathrm{KL}}(\Psi\|\Phi)=\Sc_{\mathrm{KL0}}(\Psi\|\Phi)-\int\Psi+\int
\Phi=\Sc_{\mathrm{KL0}}(\Psi\|\Phi)-\int\Psi+\frac{v^*v}{\|v^*B\|^2}\eeq
and the minimization of $\Sc_{\mathrm{KL}}(\Psi\|\Phi)$ is
equivalent to the minimization of
$\Sc_{\mathrm{KL0}}(\Psi\|\Phi)$. We conclude that the two
solutions coincide when $A$ is singular. If $A$ is instead
non-singular, the two solutions are typically different. For
instance, let
\eq \label{es_scolastico}A=\left[%
\begin{array}{cc}
  \frac{1}{2} & 0 \\
  -\sqrt 6+\sqrt\frac{8}{3} \hspace{0.2cm}& \frac{1}{3} \\
\end{array}%
\right],\;\; B=\left[%
\begin{array}{c}
  \frac{\sqrt 3}{2} \\
  \frac{\sqrt{2}}{3} \\
\end{array}%
\right],\;\;\Psi\equiv 1.\eeq It is easy to see that $\int GG^*$
is the unique solution to the Lyapunov equation $P-APA^T=BB^T$,
and in this case $P=I$. Thus $\int G\Psi G^*=I$, and $\Psi\in\Sts$
is compatible with the second-order statistics. Since
$\Sc_{\mathrm{KL}}(\Psi\|\Phi)$ is a distance measure among
spectral densities, our solution coincides with $\Psi$ (in fact
condition (\ref{condizione_der_prima_kl_raz}) holds for
$\Lambda^\circ=0$). When $\Psi\equiv 1$, the solution with
$\Sc_{KL0}(\Psi\|\Phi)$ can be expressed in the closed form
$\Phi_{\mathrm{KL0}}=(G^*B(B^*B)^{-1}B^*G)^{-1}$, see
\cite{GEORGIOU_SPECTRAL_ANALYSIS_2002}. Substituting the
parameters (\ref{es_scolastico}) in $\Phi_{\mathrm{KL0}}$, we
obtain \eq \Phi_{\mathrm{KL0}}(z)=\frac{42z^2 -245z  + 434
-245z^{-1} + 42z^{-2}}{ -175z + 370 -175z^{-1}}\eeq which is
different from the compatible prior. We conclude that our solution
preserves the approximation-feature (i.e. the solution is as close
as possible to $\Psi$) also in the case wherein $A$ is invertible.

\section{Case $1<\nu<\infty$}\label{section_THREE_aalpha} The
corresponding Lagrangian is \eqn
L(\Phi, \Lambda)&=&\Danu(\Phi\|\Psi)-\frac{\nu}{\nu-1}\int \Psi+\Sp{\int G\Phi G^*-I}{\Lambda}\nn\\
&=&\int \left[\frac{\nu^2}{1-\nu}
\Phi^{\frac{\nu-1}{\nu}}\Psi^{\frac{1}{\nu}}+\nu\Phi+G^*\Lambda
G\Phi\right]-\tr(\Lambda).\eeqn Similarly to the previous case we
can impose that $\Lambda\in\Rgamma$. Since $\Danu(\Phi\|\Psi)$ is
strictly convex with respect to $\Phi\in\Sts$, then
$L(\cdot,\Lambda)$ is strictly convex with respect to
$\Phi\in\Sts$. Accordingly, the corresponding unconstrained
minimization problem admits a unique solution given by
annihilating the first variation of $L(\cdot,\Lambda)$ in each
direction $\delta\Phi\in L_\infty(\Ts)$: \eqn
\label{variazione_delta_phi_alpha}\delta L(\Phi,\Lambda;\delta
\Phi)&=&\int\left[\left(-\nu
\Phi^{-\frac{1}{\nu}}\Psi^{\frac{1}{\nu}}+\nu+ G^*\Lambda
G\right)\delta\Phi\right]. \eeqn Note that, $-\nu
\Phi^{-\frac{1}{\nu}}\Psi^{\frac{1}{\nu}}+\nu+ G^*\Lambda G\in
L_\infty(\Ts)$. Hence, (\ref{variazione_delta_phi_alpha}) is zero
$\forall \delta\Phi\in L_\infty(\Ts)$ if and only if \eq -\nu
\Phi^{-\frac{1}{\nu}}\Psi^{\frac{1}{\nu}}+\nu+ G^*\Lambda G\equiv
0.\eeq Therefore, the unique minimum point of $L(\cdot,\Lambda)$
has the form
\eq\label{ottimo_alpha}\Phi(\Lambda):=\frac{\Psi}{\left(1+\frac{1}{\nu}
G^*\Lambda G\right)^\nu}.\eeq Since
$\Phi^{-\frac{1}{\nu}}\Psi^{\frac{1}{\nu}}\in\Sts $, the set of
admissible Lagrange multipliers  is
\eq\Lag:=\Set{\Lambda\in\Qc_n\;|\; 1+\frac{1}{\nu}G^*\Lambda G>0\;
\hbox{on} \;\Ts}\cap \Rgamma.\eeq

Also in this case $\Phi(\Lambda)$ is rational because $\Psi$ and
$G$ are rational, and $\nu\in\Zs$ is such that $\nu>1$. Moreover,
\eq\deg[\Phi(\Lambda)]\leq \deg[\Psi]+2n\nu.\eeq

It remains to be shown that the {\em dual functional}
\eq\label{funz_dual_nu}
J(\Lambda):=-L(\Phi(\Lambda),\Lambda)=\frac{\nu}{\nu-1}\int\Psi
\left(1+\frac{1}{\nu} G^*\Lambda G\right)^{1-\nu}+\tr(\Lambda)\eeq
admits a minimum point $\Lambda^\circ$ over $\Lag$. Accordingly
$\Phi(\Lambda^\circ)$ is the unique solution to Problem
\ref{problema_THREE_I}. \teo \label{teorema_J_stricly_convex}The
dual functional $J$ belongs to $\Cc^2(\Lag)$ and it is strictly
convex over $\Lag$.\eteo \proof The first variation of
$J(\Lambda)$ in direction $\delta \Lambda_1 \in\Qc_n$ is \eq
\label{gradiente_J_alpha} \delta J(\Lambda;\delta\Lambda_1)=-\int
\Psi\left(1+\frac{1}{\nu} G^*\Lambda
G\right)^{-\nu}G^*\delta\Lambda_1 G+\tr(\delta \Lambda_1).\eeq In
order to prove that $J(\Lambda)\in\Cc^1(\Lbg)$ we have to show
that $\delta J(\Lambda;\delta\Lambda_1)$, for any fixed $\delta
\Lambda_1$, is continuous in $\Lambda$. Consider a sequence
$M_n\in\Rgamma$ such that $M_n\rightarrow 0$ and define
$Q_N(z)=(1+\frac{1}{\nu}G(z)^* N G(z))^{-1}$ with $N\in\Qc_n$. We
know that $Q_{\Lambda+M_n}$ converges uniformly to $Q_\Lambda$,
see \cite[Lemma
5.2]{MATRICIAL_ALGTHM_RAMPONI_FERRANTE_PAVON_2009}. Thus, applying
elementwise the bounded convergence theorem, we obtain \eq
\lim_{n\rightarrow \infty } \int G \Psi Q_{\Lambda+M_n}^\nu G^*=
\int G \Psi Q_{\Lambda}^{\nu} G^*. \eeq Accordingly,
 $\delta J(\Lambda;\delta\Lambda)$ is continuous, i.e. $J$ belongs to $\Cc^1(\Lag)$.
  The second variation in direction
$\delta\Lambda_1,\delta\Lambda_2\in\Qc_{n}$ is \eq
\label{hessiano_J_alpha}\delta^2
J(\Lambda;\delta\Lambda_1,\delta\Lambda_2)=\int
\Psi\left(1+\frac{1}{\nu} G^*\Lambda
G\right)^{-\nu-1}G^*\delta\Lambda_1 GG^*\delta\Lambda_2 G\eeq and
the continuity of $\delta^2 J$ can be established by using the
previous argumentation. Finally, it remains to be shown that $J$
is strictly convex on the open set $\Lag$, i.e.
$\Hc_\Lambda(\delta\Lambda,\delta\Lambda):=\delta^2
J(\Lambda;\delta\Lambda,\delta\Lambda)$ is positive definite over
$\Lag$. Since the integrand in (\ref{hessiano_J_alpha}) is a
nonnegative function when $\delta\Lambda_1=\delta\Lambda_2$, we
have $\Hc_{\Lambda}(\delta\Lambda,\delta\Lambda)\geq 0$. If
$\Hc_{\Lambda}(\delta\Lambda,\delta\Lambda)=0$, then
$G^*\delta\Lambda G\equiv 0$ namely $\delta
\Lambda\in[\Rgamma]^\bot$. Since $\delta\Lambda\in\Rgamma$, it
follows that $\delta\Lambda =0$.
Thus, $\Hc_\Lambda(\delta\Lambda,\delta\Lambda)$ is positive definite. \qed\\
\teo The dual functional $J$ has a unique minimum point in
$\Lag$.\eteo \proof Since the solution of the dual problem over
$\Lag$, if it does exist, is unique, we only need to show that $J$
takes a minimum value on $\Lag$. First of all, note that $J$ is
continuous on $\Lag$, see Theorem \ref{teorema_J_stricly_convex}.
Secondly, we show that $\tr(\Lambda)$ is bounded from below on
$\Lag$. Since Problem \ref{problema_THREE_I} is feasible, there
exists $\Phi_I\in\Sts$ such that $\int G\Phi_I G^*=I$. Thus, \eqn
\tr(\Lambda)=\tr\left[\int G\Phi_I G^* \Lambda \right]=\int G^*
\Lambda G\Phi_I .  \eeqn Defining $\gamma=-\nu\int \Phi_I$, we
obtain \eq \tr(\Lambda)=\nu\int \left(1+\frac{1}{\nu}G^*\Lambda
G\right)\Phi_I +\gamma.\eeq Note that $\Phi_I$ is a coercive
spectrum, namely there exists a constant $\mu>0$ such that
$\Phi_I(\ej)\geq \mu $, $\forall \; \ej\in\Ts$. Since the integral
is a monotonic function, we get \eqn
\label{traccia_lambda_bounded_below}\tr(\Lambda)&\geq&\nu \mu\int
1+\frac{1}{\nu}G^*\Lambda G+\gamma>\gamma\eeqn where we have used
the fact that $\int 1+\frac{1}{\nu}G^*\Lambda G>0$ when
$\Lambda\in\Lag$. Thirdly, notice that $J(0)=\frac{\nu}{\nu-1}\int
\Psi$. Accordingly, we can restrict the search of a minimum point
to the set $\Kc:=\Set{\Lambda\in\Lag\;|\;J(\Lambda)\leq J(0)}$.
Finally, the existence of the solution to the dual problem follows
from the Weierstrass' theorem, since $\Kc$ is a compact set. In
order to prove that $\Kc$ is compact, it is sufficient to show
that:
\begin{enumerate}
    \item $\Lim{\Lambda\rightarrow \partial \Lag} J(\Lambda)=+\infty $;
    \item $\Lim{\|\Lambda\|\rightarrow \infty} J(\Lambda)=+\infty
    $.
\end{enumerate}
Point (1): Function $r_\Lambda(z):=1+\frac{1}{\nu}G(z)^*\Lambda
G(z)$ is rational. Observe that $\partial\Lag$ is the set of
$\Lambda\in\Rgamma$ such that $r_\Lambda(\ej)\geq 0$ on $\Ts$ and
there exists $\vartheta$ such that $r_\Lambda(\ej)$ is equal to
zero. Thus, for $\Lambda\rightarrow
\partial \Lag$ $r_\Lambda(z)^{-1}$ has at least one pole tending to the unit circle. Since $\nu\in\Zs$ and $\nu-1>0$, then
$r_\Lambda(z)^{1-\nu}$ has at least one pole (of order greater
than or equal to $\nu-1$) tending to $\Ts$. Since $\Psi$ is fixed
and coercive, then also $\Psi(z)r_\Lambda(z)^{1-\nu}$ has one pole
tending to the unit circle. Accordingly, $\int \Psi
r_\Lambda^{1-\nu}\rightarrow \infty$ as $\Lambda\rightarrow
\partial\Lag$. In view of  (\ref{traccia_lambda_bounded_below}),
we conclude that $J(\Lambda)=\frac{\nu}{\nu-1}\int
\Psi r_\Lambda^{1-\nu}+\tr(\Lambda)\rightarrow \infty$ as $\Lambda\rightarrow \partial\Lag$.\\
Point (2): Consider a sequence $\{\Lambda_k\}_{k\in\Ns}\in\Lag$,
such that \eq \Lim{k\rightarrow \infty}\|\Lambda_k\|=\infty.\eeq
Let $\Lambda^0_k=\frac{\Lambda_k}{\|\Lambda_k\|}$. Since $\Lag$ is
convex and $0\in\Lag$, if $\Lambda\in\Lag$ then
$\xi\Lambda\in\Lag$ \; $\forall \; \xi\in[0,1]$. Therefore
$\Lambda^0_k\in\Lag$ for $k$ sufficiently large. Let
$\eta:=\lim\inf \tr(\Lambda_k^0)$. In view of
(\ref{traccia_lambda_bounded_below}),
\eq\tr(\Lambda_k^0)=\frac{1}{\|\Lambda_k\|}\tr(\Lambda_k)>\frac{1}{\|\Lambda_k\|}\gamma\rightarrow
0,\eeq for $k\rightarrow \infty$, so $\eta\geq 0$. Thus, there
exists a subsequence of $\{\Lambda_k^0\}$ such that the limit of
its trace is equal to $\eta$. Moreover, this subsequence remains
on the surface of the unit ball $\partial
\Bc=\Set{\Lambda=\Lambda^T\;|\; \|\Lambda\|=1}$ which is compact.
Accordingly, it has a subsubsequence $\{\Lambda_{k_i}^0\}$
converging in $\partial \Bc$. Let $\Lambda_\infty\in\partial\Bc$
be its limit, thus $\Lim{i\rightarrow\infty}
\tr(\Lambda_{k_i}^0)=\tr(\Lambda_\infty)=\eta$. We now prove that
$\Lambda_\infty\in\Lag$. First of all, note that $\Lambda_\infty$
is the limit of a sequence in the finite dimensional linear space
$\Rgamma$, hence $\Lambda_\infty\in\Rgamma$. It remains to be
shown that $1+\frac{1}{\nu}G^*\Lambda_\infty G$ is positive
definite on $\Ts$. Consider the unnormalized sequence
$\{\Lambda_{k_i}\}\in\Lag$: We have that
$1+\frac{1}{\nu}G^*\Lambda_{k_i} G>0$ on $\Ts$ so that
$\frac{1}{\|\Lambda_{k_i}\|}+\frac{1}{\nu}G^*\Lambda_{k_i}^0 G$ is
also positive definite  on $\Ts$ for each $i$. Taking the limit
for $i\rightarrow \infty$, we get that $G^*\Lambda_\infty G$ is
positive semidefinite on $\Ts$ so that
$1+\frac{1}{\nu}G^*\Lambda_\infty G>0$ on $\Ts$. Hence,
$\Lambda_\infty\in\Lag$. Since Problem \ref{problema_THREE_I} is
feasible, there exists $\Phi_I\in\Sts$ such that $I=\int G\Phi_I
G^*$, accordingly \eq\label{eta} \eta=\tr(\Lambda_\infty)=\tr\int
G\Phi_I G^*\Lambda_\infty=\int \Phi_I G^*\Lambda_\infty G\geq \mu
\int G^*\Lambda_\infty G.\eeq Moreover, $G^*\Lambda_\infty G$ is
not identically equal to zero. In fact, if $G^*\Lambda_\infty G
\equiv0$, then $\Lambda_\infty\in[\Rgamma]^\bot$ and
$\Lambda_\infty\neq 0$ since it belongs to the surface of the unit
ball. This is a contradiction because $\Lambda_\infty\in\Rgamma$.
Thus, $G^*\Lambda_\infty G$ is not identically zero and $\eta>0$.
Finally, we have \eqn \Lim{k\rightarrow \infty}J(\Lambda_k)
&=&\Lim{k\rightarrow \infty}\frac{\nu}{\nu-1}\int
\Psi\left(1+\frac{1}{\nu}G^*\Lambda_k G\right)^{1-\nu}+\tr(\Lambda_k)\nn\\
&\geq &  \Lim{k\rightarrow
\infty}\|\Lambda_k\|\tr(\Lambda_k^0)\geq\eta\Lim{k\rightarrow
\infty} \|\Lambda_k\|=\infty. \eeqn \qed\\
\rem \label{remark_nu_neg} The optimal form (\ref{ottimo_alpha})
is also valid for $\nu\in\Zs$ such that $\nu\leq -1$. The
corresponding dual problem, however, may not have solution: The
minimum point for $J(\Lambda)$ may lie on $\partial\Lag$ since
$J(\Lambda)$ takes finite values on the boundary of $\Lag$.\erem
Now, it is worth comparing the related work in \cite{BETA}. Here,
the Beta divergence family \eq \Sc_\beta(\Phi\|\Psi)=\int
\frac{1}{\beta-1}(\Phi^\beta-\Phi\Psi^{\beta-1})-\frac{1}{\beta}(\Phi^\beta-\Psi^\beta),
\;\; \beta\in\Rs\setminus\{0,1\}\eeq has been considered. In
\cite{FLEXIBLE_ROBUST_CICHOKI_AMARI_2010}, it was shown that the
Beta divergence family can be obtained, up to a factor
$\beta^{-2}$, by the Alpha divergence family applying the
following transformation \eq \label{nonlinear_transf_a_b}\Phi
\mapsto \Phi^\beta,\;\; \Psi\mapsto \Psi^\beta,\;\; \alpha\mapsto
\frac{1}{\beta}.\eeq Conversely, the Alpha divergence family can
be obtained, up to a factor $\alpha^{-2}$, by the Beta divergence
family by the following transformation \eq
\label{nonlinear_transf_b_a}\Phi \mapsto \Phi^\alpha,\;\;
\Psi\mapsto \Psi^\alpha,\;\; \beta\mapsto \frac{1}{\alpha}.\eeq
Notice that, the above transformations are nonlinear accordingly
the assumptions on $\Psi$ could be different when the Beta
divergence is considered. In fact, taking the parametrization
$\beta=-\frac{1}{\nu}+1$ with $\nu\in\Zs$ such that
$1<\nu<\infty$, the corresponding optimal form of the spectrum
approximation problem is
\eq\label{forma_ottimo_beta}\Phi_\beta(\Lambda)=\left(\Psi^{-\frac{1}{\nu}}+\frac{1}{\nu}
G^*\Lambda G\right)^{-\nu}.\eeq Clearly, the assumption that
$\Psi$ is rational is not sufficient to guarantee that
$\Psi^{\frac{1}{\nu}}$, and thus also $\Phi_\beta(\Lambda)$, is
rational. Therefore, not only must $\Psi$ be rational, but even
$\Psi^{\frac{1}{\nu}}$ must be rational in order that
$\Phi_\beta(\Lambda)$ is rational. In this situation
$\deg[\Phi_\beta(\Lambda)]\leq
\nu(\deg[\Psi^{\frac{1}{\nu}}]+2n)$. We conclude that the solution
given by (\ref{ottimo_alpha}) is more appealing than the one given
by (\ref{forma_ottimo_beta}). Finally, note that the two solutions
coincide when $\Psi\equiv1$.

In \cite{BETA} the Beta divergence family has been extended for
the multichannel case: \eq \Sc_\beta(\Phi\|\Psi)=\tr \int
\frac{1}{\beta-1}(\Phi^\beta-\Phi\Psi^{\beta-1})-\frac{1}{\beta}(\Phi^\beta-\Psi^\beta)
\eeq where $\Phi$ and $\Psi$ are $\Rs^{m\times m}$-valued spectral
density functions which are bounded and coercive. Moreover, it was
shown that (\ref{forma_ottimo_beta}) also holds for the
multichannel case. Applying transformation
(\ref{nonlinear_transf_b_a}) we obtain the corresponding
multivariate Alpha divergence family: \eq
\Sc_\alpha(\Phi\|\Psi)=\tr \int \frac{1}{\alpha(\alpha-1)
}\Psi^{1-\alpha}\Phi^{\alpha}-\frac{1}{\alpha-1}
\Phi+\frac{1}{\alpha}\Psi,\;\; \alpha\in\Rs\setminus\{0,1\}.\eeq
By the Lieb's Theorem \cite[Theorem 1]{AUJLA_SIMPLE_LIEB_2011}, we
get that the map \eq (\Phi,\Psi)\mapsto \frac{1}{\alpha(\alpha-1)
}\tr(\Psi^{1-\alpha}\Phi^{\alpha})\eeq is convex with respect to
both $\Phi$ and $\Psi$ for $0<\alpha<1$. Accordingly,
$\Sc_\alpha(\Phi\|\Psi)$ is convex with respect to both $\Psi$ and
$\Phi$ for $0<\alpha<1$. By choosing the parametrization
$\alpha=-\frac{1}{\nu}+1$ with $\nu\in\Zs\setminus \{0,1\}$, we
obtain \eq \Sc_\nu(\Phi\|\Psi)=\tr \int \frac{\nu^2}{1-\nu
}\Phi^{\frac{\nu-1}{\nu}}\Psi^{\frac{1}{\nu}}+\nu
\Phi+\frac{\nu}{\nu-1}\Psi.\eeq The corresponding Lagragian for
the multivariate version of Problem \ref{problema_THREE_I} is \eq
L(\Phi, \Lambda)=\Danu(\Phi\|\Psi)-\frac{\nu}{\nu-1}\tr \int
\Psi+\Sp{\int G\Phi G^*-I}{\Lambda}.\eeq The minimum point for
$L(\cdot,\Lambda)$ must annihilate its first variation (see
\cite[Appendix]{BETA} for the differential of the map $X\mapsto
X^c$) in each direction $\delta \Phi\in\L_\infty^{m\times m}
(\Ts)$  \eqn && \delta L(\Phi,\Lambda;\delta\Phi)\nn\\ &&
\hspace{0.5cm}=\tr
\int-\nu\Psi^{\frac{1}{\nu}}\int_0^1\int_0^\infty \left(
\Phi^{\frac{\nu-1}{\nu}(1-\tau)}(\Phi+tI)^{-1}\delta\Phi(\Phi+tI)^{-1}\Phi^{\frac{\nu-1}{\nu}\tau}
\right)\de t\de\tau \nn \\ && \hspace{0.8cm} +\nu\delta\Phi
+G^*\Lambda G\delta\Phi\eeqn which implies the following condition
\eq \int_0^1\int_0^\infty
(\Phi+tI)^{-1}\Phi^{\frac{\nu-1}{\nu}\tau}\Psi^{\frac{1}{\nu}}\Phi^{\frac{\nu-1}{\nu}(1-\tau)}(\Phi+tI)^{-1}
\de t\de \tau=I+\frac{1}{\nu}G^*\Lambda G.\eeq Hence, it is not
possible to derive an optimal form for the multichannel case and
(\ref{ottimo_alpha}) is only valid for the scalar case. We
conclude that the Beta divergence family is more appealing than
the Alpha divergence family when multivariate processes are
considered.

\section{Case
$\nu=\infty$}\label{section_THREE_KLinf} In this case we minimize
the {\em Kullback-Leibler} divergence with respect to the first
argument. Thus, the corresponding solution follows the principle
of minimum discrimination information. Its {\em Lagrangian} is
\eqn
L(\Phi,\Lambda)&:=&\Sc_{\mathrm{KL}}(\Phi\|\Psi)-\int\Psi+\Sp{\int
G\Phi G^*-I}{\Lambda}\nn\\ &=&  \int \Phi(\log
(\Phi)-\log(\Psi))-\Phi+G^*\Lambda G\Phi-\tr(\Lambda)\nn \eeqn
where $\Lambda\in\Rgamma$ as in the previous cases. Since
$L(\cdot,\Lambda)$ is strictly convex over $\Sts$, its unique
minimum point is given by annihilating its first directional
derivative for each $\delta\Phi\in \Lm$: \eq \delta
L(\Phi,\Lambda;\delta \Phi)=\int[ \log(\Phi)-\log(\Psi)+G^*\Lambda
G]\delta \Phi.\eeq Since $\log(\Phi)-\log(\Psi)+G^*\Lambda G\in
L_\infty(\Ts)$, the minimum point for $L(\cdot,\Lambda)$ is
\eq\label{Phi_Lambda_inf} \Phi(\Lambda):=\e{\log(\Psi)-G^*\Lambda
G}\eeq and the set of admissible Lagrange multipliers is
$\Rgamma$. Note that (\ref{Phi_Lambda_inf}) is an exponential
solution. Thus, it is not rational though $\Psi$ is rational.
Finally, we show that the {\em dual functional} \eq
\label{funz_dual_KL_inf}J(\Lambda)=-L(\Phi(\Lambda),\Lambda)=\int
\Psi\e{-G^*\Lambda G}+\tr(\Lambda)\eeq admits a minimum point
$\Lambda^\circ$ over $\Rgamma$. Hence, $\Phi(\Lambda^\circ)$ is
the unique solution to Problem \ref{problema_THREE_I}. \teo The
dual functional $J$ belongs to $\Cc^2(\Rgamma)$ and it is strictly
convex over $\Rgamma$.\eteo \proof The first and the second
variation of $J(\Lambda)$ in direction $\delta \Lambda \in\Qc_n$
are  \eqn \delta J(\Lambda;\delta\Lambda)&=&-\int \Psi\e{-
G^*\Lambda G}G^*\delta\Lambda G+\tr(\delta \Lambda)\nn\\ \delta^2
J(\Lambda;\delta\Lambda,\delta\Lambda)&=&\int \Psi\e{-G^*\Lambda
G}(G^*\delta\Lambda G)^2,\eeqn respectively. Similarly to the
previous cases, in order to prove that
$J(\Lambda)\in\Cc^1(\Rgamma)$ we consider a sequence
$M_n\in\Rgamma$ such that $M_n\rightarrow 0$ and define
$Q_N(z)=\e{-G(z)^*NG(z)}$ with $N\in\Qc_n$. By using similar
argumentations in the proof of Lemma 5.2 in
\cite{MATRICIAL_ALGTHM_RAMPONI_FERRANTE_PAVON_2009}, it is
possible to prove that  $Q_{\Lambda+M_n}$ converges uniformly to
$Q_{\Lambda}$. Thus, applying the bounded convergence theorem, we
obtain \eq \lim_{n\rightarrow \infty } \int G \Psi Q_{\Lambda+M_n}
G^*= \int G \Psi Q_{\Lambda} G^*. \eeq Accordingly
 $\delta J(\Lambda ;\delta\Lambda)$, once fixed $\delta\Lambda$, is continuous in $\Lambda$, i.e. $J$ belongs to $\Cc^1(\Rgamma)$.
In similar way we can establish the continuity of $\delta^2
J(\Lambda;\delta\Lambda,\delta\Lambda)$. Finally, note that
$\Hc_\Lambda(\delta\Lambda,\delta\Lambda):=\delta^2
J(\Lambda;\delta\Lambda,\delta\Lambda)\geq 0$ and
$\Psi\e{-G^*\Lambda G}\in\Sts$. Thus, $\delta^2
J(\Lambda;\delta\Lambda,\delta\Lambda)= 0$ implies that
$G^*\delta\Lambda G\equiv 0$, i.e.
$\delta\Lambda\in[\Rgamma]^\bot$. Since $\delta\Lambda\in\Rgamma$,
we get $\delta\Lambda=0$. We conclude that the {\em Hessian} of
$J$ is positive definite, and thus $J$ is strictly convex over
$\Rgamma$. \qed\\ \teo \label{teorema_esistenza_KL_inf} $J$ admits
a unique minimum point over $\Rgamma$. \eteo  \proof Also in this
case it is sufficient to show that $J$ takes a minimum value on
$\Rgamma$. Firstly, note that $J(0)=\int \Psi$. Hence, the search
of the minimum point over $\Rgamma$ is equivalent to minimize $J$
over the closed set $\Kc:=\{\Lambda\in\Rgamma\;|\; J(\Lambda)\leq
J(0)\}$. We want to show that $\Kc$ is bounded and accordingly
compact. Then, by Weierstrass' theorem we conclude that $J$ admits
a minimum point over $\Kc$. To show that $\Kc$ is bounded, we
prove that
\eq\lim_{\|\Lambda\|\rightarrow\infty}J(\Lambda)=\infty.\eeq To
this aim, note that Problem \ref{problema_THREE_I} is feasible.
Thus, there exists $\Phi_I\in\Sts$ such that $\int G\Phi_I G^*=I$.
Accordingly, \eqn \tr(\Lambda)=\int \Phi_IG^*\Lambda G.\eeqn Let
us consider a sequence $\{\Lambda_k\}_{k\in\Ns}$, $\Lambda_k
\in\Rgamma$ such that $\|\Lambda_k\|\rightarrow\infty$. Consider
the sequence $\Lambda^0_k=\frac{\Lambda_k}{\|\Lambda_k\|}$ which
is contained in the closed ball $\{\Lambda=\Lambda^T\;|\;
\|\Lambda\|=1\}$. Let $\eta=\lim\inf \tr(\Lambda_k^0)$. Note that
$|\eta|<\infty$. Consider a subsequence of $\{\Lambda_k^0\}$ such
that the limit of its trace is equal to $\eta$. Since this
subsequence is contained in a compact set, there exists a
subsubsequence $\{\Lambda^0_{k_i}\}$ having limit
$\Lambda_\infty\in\Rgamma$ with $\|\Lambda_\infty\|=1$. Clearly,
$\tr(\Lambda_\infty)=\eta$. Note that $\Lambda_\infty$ is not
equal to the null matrix because $\|\Lambda_\infty\|=1$. Moreover
$G^*\Lambda_\infty G\neq 0$ because $\Lambda_\infty\in\Rgamma$. If
$G^*\Lambda_\infty G\geq 0$, then \eq\eta=\tr(\Lambda_\infty)=\int
\Phi_IG^*\Lambda_\infty G\geq \mu\int G^*\Lambda_\infty G>0,\eeq
where we have exploited the fact that $\Phi_I\in\Sts$, i.e. there
exists $\mu>0$ such that $\Phi_I(\ej)\geq \mu$ $\forall \;
\e{j\vartheta}\in\Ts$. Accordingly, \eq\lim_{k\rightarrow\infty}
J(\Lambda_k)\geq \lim_{k\rightarrow\infty} \|\Lambda_k\|\int
\Phi_I G^*\Lambda_k^0G\geq \eta\lim_{k\rightarrow \infty}
\|\Lambda_k\|=\infty.\eeq In the remaining possible case, there
exists $\bar\vartheta$ such that $G^*(\bar\vartheta)\Lambda_\infty
G(\bar\vartheta)<0$. Thus, $G^*(\bar\vartheta)\Lambda_k
G(\bar\vartheta)\rightarrow -\infty$, and accordingly $\int
\Psi\e{-G^*\Lambda_k G}\rightarrow \infty$, as $k\rightarrow
\infty$. Moreover the latter dominates the term
$\tr(\Lambda_k)=\int\Phi_I G^*\Lambda_kG $. Accordingly, we
conclude that $J(\Lambda_k)\rightarrow \infty$ as $k\rightarrow
\infty$.
 \qed\\ Note that Theorem \ref{teorema_esistenza_KL_inf} can also be proven by using homotopy-like methods as
in \cite{RELATIVE_ENTROPY_GEORGIOU_2006}.

\section{Features of the family of solutions}
Let $\Phi_\nu(\Lambda)$ denote the optimal form
(\ref{ottimo_alpha}) and $\Phi_\infty(\Lambda)$ denote the optimal
form (\ref{Phi_Lambda_inf}).  In this Section, we want to show
$\Phi_\nu(\Lambda)\rightarrow \Phi_\infty(\Lambda)$ uniformly on
$\Ts$ as $\nu\rightarrow \infty$. By exploiting the limit
$\Lim{\nu\rightarrow \infty} (1+\frac{1}{\nu}x)^\nu=\e{x}$, we
obtain the following pointwise limit \eqn \Lim{\nu\rightarrow
\infty}\Phi_{\nu}(\Lambda)=\Lim{\nu\rightarrow \infty}
\frac{\Psi}{(1+\frac{1}{\nu}G^*\Lambda G)^\nu}=\Psi \e{-G^*\Lambda
G}=\Phi_\infty(\Lambda).\eeqn \prop \label{corollario_conv_kl}
Assume that $\|\Lambda\|<\infty$, then $\Phi_\nu(\Lambda)$
converges uniformly to $\Phi_\infty(\Lambda)$ on $\Ts$ as
$\nu\rightarrow \infty$. \eprop \proof Since $\Psi\in\Sts$, there
exists a constant $C_\Psi$ such that $\Psi(\e{i\vartheta})\leq
C_\Psi$ on $\Ts$. Let $\nu > \max\{1,M_G\}$ and $
M_G:=\max_{\theta} \|G^*(\e{i\vartheta})\|\|\Lambda\|
\|G(\e{i\vartheta})\|$. Then \eqn \sup_\vartheta
\left|\Phi_\nu(\Lambda)-\Phi_\infty(\Lambda)\right|&=&
\sup_\vartheta |\Psi|\left|(1+\frac{1}{\nu}G^*\Lambda
G)^{-\nu}-\e{-G^*\Lambda G}\right|\nn\\ &\leq & C_\Psi
\sup_\vartheta \left|(1+\frac{1}{\nu}G^*\Lambda
G)^{-\nu}-\e{-G^*\Lambda G}\right|.\eeqn Let us consider a first
order Taylor expansion of $(1+\frac{1}{\nu}G^*\Lambda G)^{-\nu}$:
\eq (1+\frac{1}{\nu}G^*\Lambda G)^{-\nu}=\e{-G^*\Lambda
G}+f(\xi,\vartheta)\frac{1}{\nu}\eeq for a certain $0 \leq \xi\leq
\frac{1}{\nu}$.  Here,  \eqn
f(\xi,\vartheta)&:=&\frac{\partial}{\partial \xi} (1+\xi
G^*\Lambda G)^{-\frac{1}{\xi}}\nn\\&=&(1+\xi G^*\Lambda
G)^{-\frac{1}{\xi}} \left(-\frac{G^*\Lambda G}{\xi(1+\xi
G^*\Lambda G)}+\frac{\log(1+\xi G^*\Lambda G)}{\xi^2}\right)\eeqn
when $0 < \xi\leq \frac{1}{\nu}$ and we extend it by continuity in
$\xi=0$ \eq f(0,\vartheta):=\lim_{\xi\rightarrow
0}f(\xi,\vartheta)=\e{-G^*\Lambda G} (G^*\Lambda G)^2.\eeq
Accordingly,\eqn \sup_\vartheta
\left|\Phi_\nu(\Lambda)-\Phi_\infty(\Lambda)\right|&\leq&
\frac{C_\Psi}{\nu} \sup_\vartheta |f(\xi,\vartheta)|.\eeqn Since
$\xi\leq \frac{1}{\nu}< \frac{1}{M_G}$, $1+\xi
G^*(\e{i\vartheta})\Lambda G(\e{i\vartheta})>0$ on  $\Ts$.
Accordingly, $f(\xi,\vartheta)$ is continuous over the compact set
$[0,\frac{1}{\nu}]\times [0,2\pi]$, and by Weierstrass' theorem it
admits minimum and maximum over such a set:
\eqn K_1(\nu)&=&\max_{\xi\in[0,\frac{1}{\nu}],\vartheta\in[0,2\pi]} f(\xi,\vartheta)\nn\\
K_2(\nu)&=&\min_{\xi\in[0,\frac{1}{\nu}],\vartheta\in[0,2\pi]}
f(\xi,\vartheta).\eeqn Hence, \eq \sup_\vartheta
\left|\Phi_\nu(\Lambda)-\Phi_\infty(\Lambda)\right| \leq
\frac{C_\Psi}{\nu}K\eeq where $K:=\max\{|K_1|,|K_2|\}$. We
conclude that $\Phi_\nu(\Lambda)\rightarrow \Phi_\infty(\Lambda)$
uniformly
on $\Ts$. \qed\\

Note that the optimal solution to the dual problem changes by
changing $\nu$. Let $\Lambda^\circ_\nu$ and $\Lambda^\circ_\infty$
be the optimal Lagrange multipliers of (\ref{funz_dual_nu}) and
(\ref{funz_dual_KL_inf}), respectively. By Proposition
\ref{corollario_conv_kl} we cannot conclude that
$\Phi_\nu(\Lambda_\nu^\circ)\rightarrow\Phi_\infty(\Lambda_\infty^\circ)$
uniformly on $\Ts$ as $\nu\rightarrow \infty$. However,
simulations suggest this conjecture is true. To illustrate this
fact, we analyze the case of the ARMA process considered in
\cite[Section VIIB]{MATRICIAL_ALGTHM_RAMPONI_FERRANTE_PAVON_2009}
with spectral density \eq \Omega(z)=\frac{ z^5
+1.1z^4+0.08z^3-0.15z^4}{z^5-0.5z^4+0.42z^3-0.602z^2 + 0.0425
z-0.1192}.\eeq  We
choose as filters bank \eq G(z)=\left[%
\begin{array}{ccc}
  z^{-6} & \ldots & z^{-1} \\
\end{array}%
\right]^T\eeq and the corresponding output covariance is
\eq\Sigma=\int G\Omega G^*\simeq\left[
\begin{array}{cccccc}
    5.58 \;\; & \;\;3.74 \; \;& \;\;1.85 \;\; &  \;\;   2.63\;\; & \;\; 3 \;\;&\; \;2.01\\
    3.74   &   5.58   &   3.74  &    1.85  &    2.63  &    3\\
    1.85   &   3.74   &   5.58  &    3.74   &   1.85   &   2.63\\
    2.63   &   1.85   &   3.74   &   5.58   &   3.74   &   1.85\\
    3   &   2.63  &   1.85   &   3.74   &   5.58  &    3.74\\
    2.01   &   3   &   2.63   &   1.85   &   3.74  &    5.58\\
\end{array}%
\right]. \eeq The {\em a priori} information on the ARMA process
is given by the {\em prior} \eq \Psi(z)=\frac{z}{z-0.82}.\eeq
\begin{figure}[htbp]
\begin{center}
\includegraphics[width=12cm]{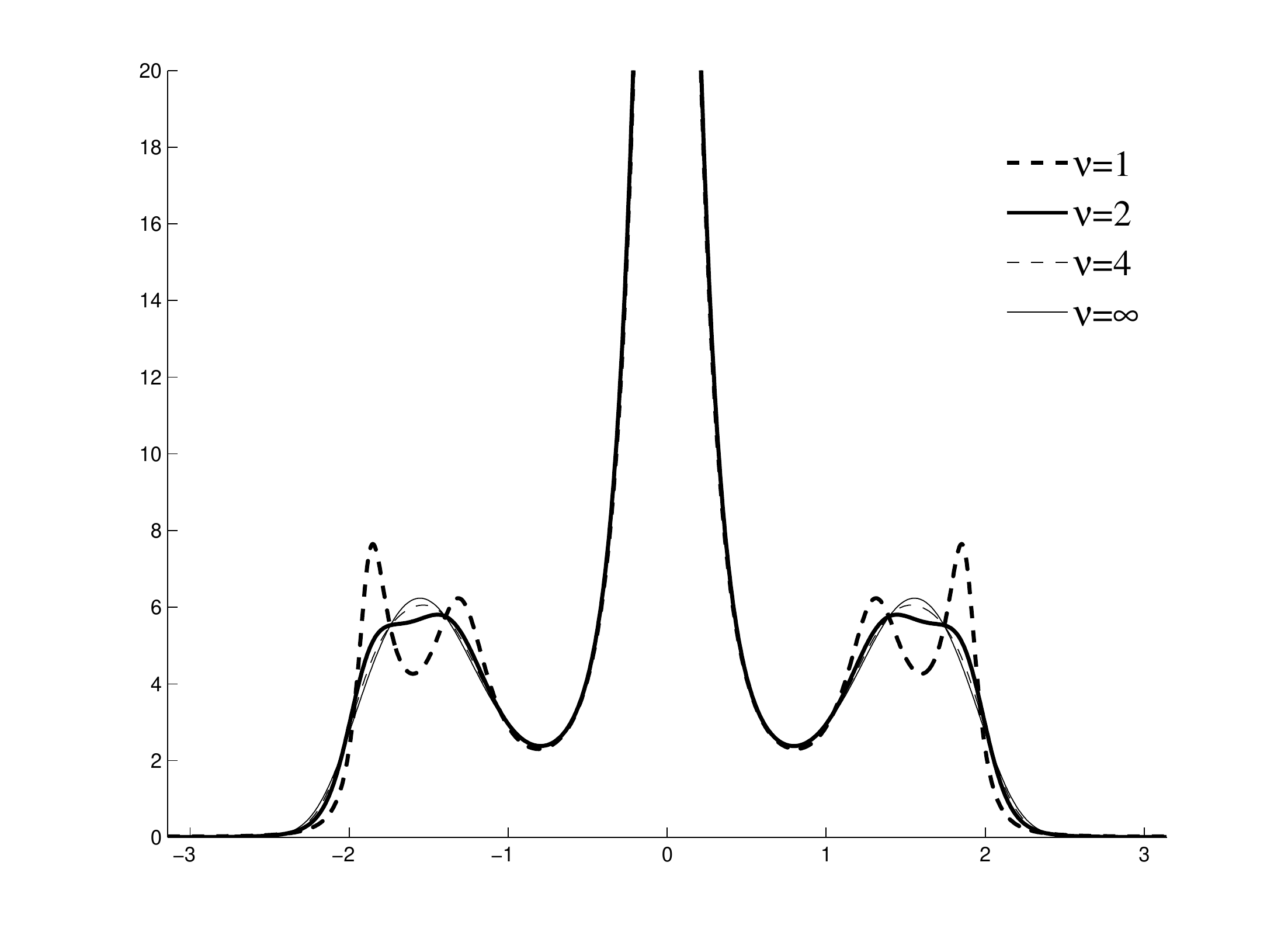}
\end{center}
\vspace{-6mm} \caption{Approximations of $\Psi$ with different
values of $\nu$.}\label{foto2}
\end{figure}
In Figure \ref{foto2} the approximations of $\Psi$ compatible with
$\Sigma$ for $\nu=1$, $\nu=2$, $\nu=4$ and $\nu=\infty$ are
depicted. It is clear that increasing $\nu$ the corresponding
solution approaches the MinxEnt solution ($\nu=\infty$), and for
$\nu=4$ it is pretty similar to the MinxEnt one. We have conducted
other numerical experiments and we observed the same behaviour as
$\nu$ increases.

\section{Conclusion}
We analyzed a spectrum approximation problem based on a suitable
parametrization of the Alpha divergence family. Here, we make the
mild assumption that the prior $\Psi$ is rational. When the
parameter $\nu\in\Zs$ is such that $1\leq \nu<\infty$, the
corresponding family of solutions is rational with an upper bound
on the degree equal to $\deg[\Psi]+2n\nu$. Moreover, the solution
having the smallest upper bound is given by minimizing the {\em
Kullback-Leibler} divergence with respect to the second argument
(case $\nu=1$). Such solution also generalizes the approximation
presented in \cite{KL_APPROX_GEORGIUO_LINDQUIST} which only holds
when the matrix $A$ is singular. When  $\nu$ tends to infinity the
solution of this family uniformly converges on $\Ts$ to an
``exponential-type'' solution having the same structure as the
{\em minimum discrimination information} solution (MinxEnt)
obtained with $\nu=\infty$. Moreover, numerical experiments show
that solutions with $\nu$ large are almost equal to the MinxEnt
solution. Hence, the family of solutions based on the Alpha
divergence yields a concrete way to approximate the MinxEnt
solution with a rational solution.

\section*{Acknowledgements}
This work was supported by University of Padova under the project
``A Unifying Framework for Spectral Estimation and Matrix
Completion: A New Paradigm for Identification, Estimation, and
Signal Processing".
 
%\bibliography{biblio}

\end{document}